\documentclass[12pt]{article}

\usepackage{amsmath,amssymb,amsthm,amscd,a4wide,upref}

\usepackage[enableskew]{youngtab}

\usepackage{hyperref}
\hypersetup{colorlinks,citecolor=blue,filecolor=black,linkcolor=blue,urlcolor=blue}
\usepackage{enumerate}
\usepackage{mathtools}
\usepackage{enumitem}

\usepackage{lmodern}     
\usepackage[T1]{fontenc}

\makeatletter
\DeclareTextCommand{\textprime}{\encodingdefault}{%
  \mbox{$\m@th'\kern-\scriptspace$}%
}
\makeatother

\begin{document}


\newcommand{\ad}{{\rm ad}}
\newcommand{\cri}{{\rm cri}}
\newcommand{\ext}{{\rm ext}}
\newcommand{\row}{{\rm row}}
\newcommand{\sh}{{\rm{sh}}}
\newcommand{\col}{{\rm col}}
\newcommand{\End}{{\rm{End}\ts}}
\newcommand{\Rep}{{\rm{Rep}\ts}}
\newcommand{\Hom}{{\rm{Hom}}}
\newcommand{\Mat}{{\rm{Mat}}}
\newcommand{\ch}{{\rm{ch}\ts}}
\newcommand{\chara}{{\rm{char}\ts}}
\newcommand{\diag}{{\rm diag}}
\newcommand{\st}{{\rm st}}
\newcommand{\non}{\nonumber}
\newcommand{\wt}{\widetilde}
\newcommand{\wh}{\widehat}
\newcommand{\ol}{\overline}
\newcommand{\ot}{\otimes}
\newcommand{\la}{\lambda}
\newcommand{\La}{\Lambda}
\newcommand{\De}{\Delta}
\newcommand{\al}{\alpha}
\newcommand{\be}{\beta}
\newcommand{\ga}{\gamma}
\newcommand{\Ga}{\Gamma}
\newcommand{\ep}{\epsilon}
\newcommand{\ka}{\kappa}
\newcommand{\vk}{\varkappa}
\newcommand{\vt}{\vartheta}
\newcommand{\si}{\sigma}
\newcommand{\vs}{\varsigma}
\newcommand{\vp}{\varphi}
\newcommand{\de}{\delta}
\newcommand{\ze}{\zeta}
\newcommand{\om}{\omega}
\newcommand{\Om}{\Omega}
\newcommand{\ee}{\epsilon^{}}
\newcommand{\su}{s^{}}
\newcommand{\hra}{\hookrightarrow}
\newcommand{\ve}{\varepsilon}
\newcommand{\ts}{\,}
\newcommand{\vac}{\mathbf{1}}
\newcommand{\vacu}{|0\rangle}
\newcommand{\di}{\partial}
\newcommand{\qin}{q^{-1}}
\newcommand{\tss}{\hspace{1pt}}
\newcommand{\Sr}{ {\rm S}}
\newcommand{\U}{ {\rm U}}
\newcommand{\BL}{ {\overline L}}
\newcommand{\BE}{ {\overline E}}
\newcommand{\BP}{ {\overline P}}
\newcommand{\cR}{\check{R}}
\newcommand{\AAb}{\mathbb{A}\tss}
\newcommand{\CC}{\mathbb{C}\tss}
\newcommand{\KK}{\mathbb{K}\tss}
\newcommand{\QQ}{\mathbb{Q}\tss}
\newcommand{\SSb}{\mathbb{S}}
\newcommand{\TT}{\mathbb{T}\tss}
\newcommand{\ZZ}{\mathbb{Z}\tss}
\newcommand{\DY}{ {\rm DY}}
\newcommand{\X}{ {\rm X}}
\newcommand{\Y}{ {\rm Y}}
\newcommand{\Z}{{\rm Z}}
\newcommand{\Ac}{\mathcal{A}}
\newcommand{\Dc}{\mathcal{D}}
\newcommand{\Lc}{\mathcal{L}}
\newcommand{\Mc}{\mathcal{M}}
\newcommand{\Pc}{\mathcal{P}}
\newcommand{\Qc}{\mathcal{Q}}
\newcommand{\Rc}{\mathcal{R}}
\newcommand{\Sc}{\mathcal{S}}
\newcommand{\Tc}{\mathcal{T}}
\newcommand{\Bc}{\mathcal{B}}
\newcommand{\Cc}{\mathcal{C}}
\newcommand{\Ec}{\mathcal{E}}
\newcommand{\Fc}{\mathcal{F}}
\newcommand{\Gc}{\mathcal{G}}
\newcommand{\Hc}{\mathcal{H}}
\newcommand{\Uc}{\mathcal{U}}
\newcommand{\Vc}{\mathcal{V}}
\newcommand{\Wc}{\mathcal{W}}
\newcommand{\Xc}{\mathcal{X}}
\newcommand{\Yc}{\mathcal{Y}}
\newcommand{\Ar}{{\rm A}}
\newcommand{\Br}{{\rm B}}
\newcommand{\Ir}{{\rm I}}
\newcommand{\Fr}{{\rm F}}
\newcommand{\Jr}{{\rm J}}
\newcommand{\Or}{{\rm O}}
\newcommand{\GL}{{\rm GL}}
\newcommand{\Spr}{{\rm Sp}}
\newcommand{\Rr}{{\rm R}}
\newcommand{\Zr}{{\rm Z}}
\newcommand{\gl}{\mathfrak{gl}}
\newcommand{\middd}{{\rm mid}}
\newcommand{\ev}{{\rm ev}}
\newcommand{\Pf}{{\rm Pf}}
\newcommand{\Norm}{{\rm Norm\tss}}
\newcommand{\oa}{\mathfrak{o}}
\newcommand{\spa}{\mathfrak{sp}}
\newcommand{\osp}{\mathfrak{osp}}
\newcommand{\f}{\mathfrak{f}}
\newcommand{\g}{\mathfrak{g}}
\newcommand{\h}{\mathfrak h}
\newcommand{\n}{\mathfrak n}
\newcommand{\m}{\mathfrak m}
\newcommand{\z}{\mathfrak{z}}
\newcommand{\Zgot}{\mathfrak{Z}}
\newcommand{\p}{\mathfrak{p}}
\newcommand{\sll}{\mathfrak{sl}}
\newcommand{\whg}{\wh{\g}}
\newcommand{\gll}{\g^{}_{\ell}}
\newcommand{\gllh}{\wh{\g}^{}_{\ell}}
\newcommand{\gllm}{\g^{}_{\ell,\ell'}}
\newcommand{\glls}{\g^*_{\ell}}
\newcommand{\agot}{\mathfrak{a}}
\newcommand{\bgot}{\mathfrak{b}}
\newcommand{\qdet}{ {\rm qdet}\ts}
\newcommand{\tra}{ {\rm t}}
\newcommand{\Ber}{ {\rm Ber}\ts}
\newcommand{\HC}{ {\mathcal HC}}
\newcommand{\cdet}{{\rm cdet}}
\newcommand{\rdet}{{\rm rdet}}
\newcommand{\tr}{ {\rm tr}}
\newcommand{\gr}{ {\rm gr}\ts}
\newcommand{\str}{ {\rm str}}
\newcommand{\loc}{{\rm loc}}
\newcommand{\Fun}{{\rm{Fun}\ts}}
\newcommand{\Gr}{{\rm G}}
\newcommand{\sgn}{ {\rm sgn}\ts}
\newcommand{\sign}{{\rm sgn}}
\newcommand{\ba}{\bar{a}}
\newcommand{\bb}{\bar{b}}
\newcommand{\bi}{\bar{\imath}}
\newcommand{\bj}{\bar{\jmath}}
\newcommand{\bk}{\bar{k}}
\newcommand{\bl}{\bar{l}}
\newcommand{\hb}{\mathbf{h}}
\newcommand{\Sym}{\mathfrak S}
\newcommand{\fand}{\quad\text{and}\quad}
\newcommand{\Fand}{\qquad\text{and}\qquad}
\newcommand{\For}{\qquad\text{or}\qquad}
\newcommand{\for}{\quad\text{or}\quad}
\newcommand{\grpr}{{\rm gr}^{\tss\prime}\ts}
\newcommand{\degpr}{{\rm deg}^{\tss\prime}\tss}
\newcommand{\bideg}{{\rm bideg}\ts}

\renewcommand{\theequation}{\arabic{section}.\arabic{equation}}

\numberwithin{equation}{section}

\newtheorem{thm}{Theorem}[section]
\newtheorem{lem}[thm]{Lemma}
\newtheorem{prop}[thm]{Proposition}
\newtheorem{cor}[thm]{Corollary}
\newtheorem{conj}[thm]{Conjecture}
\newtheorem*{mthm}{Main Theorem}
\newtheorem*{mthma}{Theorem A}
\newtheorem*{mthmb}{Theorem B}
\newtheorem*{mthmc}{Theorem C}
\newtheorem*{mthmd}{Theorem D}

\theoremstyle{definition}
\newtheorem{defin}[thm]{Definition}

\theoremstyle{remark}
\newtheorem{remark}[thm]{Remark}
\newtheorem{example}[thm]{Example}
\newtheorem{examples}[thm]{Examples}

\newcommand{\bth}{\begin{thm}}
\renewcommand{\eth}{\end{thm}}
\newcommand{\bpr}{\begin{prop}}
\newcommand{\epr}{\end{prop}}
\newcommand{\ble}{\begin{lem}}
\newcommand{\ele}{\end{lem}}
\newcommand{\bco}{\begin{cor}}
\newcommand{\eco}{\end{cor}}
\newcommand{\bde}{\begin{defin}}
\newcommand{\ede}{\end{defin}}
\newcommand{\bex}{\begin{example}}
\newcommand{\eex}{\end{example}}
\newcommand{\bes}{\begin{examples}}
\newcommand{\ees}{\end{examples}}
\newcommand{\bre}{\begin{remark}}
\newcommand{\ere}{\end{remark}}
\newcommand{\bcj}{\begin{conj}}
\newcommand{\ecj}{\end{conj}}

\newcommand{\bal}{\begin{aligned}}
\newcommand{\eal}{\end{aligned}}
\newcommand{\beq}{\begin{equation}}
\newcommand{\eeq}{\end{equation}}
\newcommand{\ben}{\begin{equation*}}
\newcommand{\een}{\end{equation*}}

\newcommand{\bpf}{\begin{proof}}
\newcommand{\epf}{\end{proof}}

\def\beql#1{\begin{equation}\label{#1}}

\newcommand{\Res}{\mathop{\mathrm{Res}}}

\title{\Large\bf The $q$-immanants and higher quantum Capelli identities}

\author{{Naihuan Jing,\quad Ming Liu\quad and\quad Alexander Molev}}

\date{} 
\maketitle


\begin{abstract}
We construct polynomials $\SSb_{\mu}(z)$ parameterized
by Young diagrams $\mu$, whose coefficients are
central elements of the quantized enveloping algebra $\U_q(\gl_n)$. Their constant terms
coincide with the central elements provided by the general construction of Drinfeld and Reshetikhin.
For another special value of $z$, we get $q$-analogues of
Okounkov's quantum immanants for $\gl_n$. We show that the Harish-Chandra image of $\SSb_{\mu}(z)$
is a factorial Schur polynomial. We also prove quantum analogues of the higher
Capelli identities and derive Newton-type identities.


\end{abstract}

\section{Introduction}
\label{sec:int}

The {\em quantum $\mu$-immanants} $\SSb_{\mu}$ are elements of the center of the universal enveloping
algebra $\U(\gl_n)$ parameterized by Young diagrams $\mu$ with at most $n$ rows.
They were introduced by
Okounkov~\cite{o:qi} by
relying on the Schur--Weyl duality between $\gl_n$ and the symmetric group $\Sym_m$.
He showed that the elements $\SSb_{\mu}$ form a basis of the center.
Moreover, the images of $\SSb_{\mu}$ under the Harish-Chandra isomorphism are
the {\em shifted Schur polynomials} $s^*_{\mu}$ which were studied in detail in the work by
Okounkov and Olshanski~\cite{oo:ss}. They essentially coincide
with a particular family of the {\em factorial Schur polynomials}
which form a distinguished basis of the algebra
of symmetric polynomials in $n$ variables as introduced in \cite{gg:nt} and \cite{m:sft}.

It was also shown in \cite{o:qi} that the images
of the quantum $\mu$-immanants in the algebra of polynomial coefficient differential operators
admit a simple factorized form thus generalizing the celebrated Capelli identity~\cite{c:zc}.
It was pointed out in \cite{jlz:gc}, that these {\em higher Capelli identities} were previously
obtained by Williamson~\cite{w:so} in a different form. The paper \cite{jlz:gc} also gives generalized
Capelli-type identities, containing many previous versions as particular cases; see {\em loc.~cit.}
for a review.

Our goal in this paper is to introduce $q$-analogues of the quantum immanants as elements
of the center of the quantized enveloping algebra $\U_q(\gl_n)$ and prove a quantum version
of the higher Capelli identities.

It is well-known that
the center of $\U_q(\gl_n)$ is isomorphic to the subalgebra of the algebra of Laurent
polynomials in the variables $q^{\ell_1},\dots,q^{\ell_n}$ generated by
the algebra of symmetric polynomials $\CC[q^{2\ell_1},\dots,q^{2\ell_n}]^{\Sym_n}$
and the element $q^{-(\ell_1+\dots+\ell_n)}$ via a quantum version
of the Harish-Chandra isomorphism; see \cite{b:cq}, \cite{jl:lf}, \cite{r:af} and \cite{t:hi}.
In terms of the $RLL$-presentation~\cite{rtf:ql},
explicit generators of the center are found as
the {\em quantum traces} $\tr_q\tss L^m$ of the powers of
the generator matrix $L$ or as
the coefficients of the
{\em quantum determinant} associated with $L$.
These coefficients were used by Noumi, Umeda and Wakayama~\cite{nuw:qa} to produce $q$-analogues
of the classical Capelli identity.

In the $RLL$-presentation, the algebra $\U_q(\gl_n)$ is generated by the entries of two
triangular matrices
\ben
L^+=\begin{bmatrix}l^{\ts+}_{11}&l^{\ts+}_{12}&\dots&l^{\ts+}_{1n}\\
0&l^{\ts+}_{22}&\dots&l^{\ts+}_{2n}\\
                                   \vdots&\vdots&\ddots&  \vdots \\
                             0&0&\dots&l^{\ts+}_{nn}
\end{bmatrix}
\Fand
L^-=\begin{bmatrix}l^{\ts-}_{11}&0&\dots&0\\
l^{\ts-}_{21}&l^{\ts-}_{22}&\dots&0\\
                                   \vdots&\vdots&\ddots&  \vdots \\
                             l^{\ts-}_{n1}&l^{\ts-}_{n2}&\dots&l^{\ts-}_{nn}
                \end{bmatrix}
\een
subject to the defining relations
\begin{align}
l^+_{ii}\ts l^-_{ii}&=l^-_{ii}\ts l^+_{ii}=1,\qquad 1\leqslant i\leqslant n,
\non\\[0.3em]
R\ts L^{\pm}_1L^{\pm}_2&=L^{\pm}_2L^{\pm}_1R,\qquad
R\ts L^+_1L^-_2=L^-_2 L^+_1R,
\label{rtt}
\end{align}
where the $R$-matrix $R$ is defined by \eqref{R}
and the standard subscript notation is explained in Sec.~\ref{sec:dch}.
It follows from the defining relations that the matrix
$L=L^+(L^-)^{-1}$ satisfies the {\em reflection equation}
\ben
RL_1R_{21}L_2=L_2RL_1R_{21}
\een
which can also be written as
\beql{re}
\cR\tss L_1\cR\tss L_1=L_1\cR\tss L_1\cR
\eeq
with the use of the $R$-matrix $\cR=PR$, where $P$ is the permutation operator
in $\CC^n\ot\CC^n$ defined in \eqref{p}. We denote by $\U^{\circ}_q(\gl_n)$ the subalgebra
of $\U_q(\gl_n)$ generated by the entries of the matrix $L$.
It turns out that the abstract algebra ({\em reflection equation algebra})
generated by the entries of the matrix $L$ subject
to the defining relations \eqref{re} is `almost isomorphic' to $\U_q(\gl_n)$; see e.g.
\cite[Sec.~1.4]{jw:cr}
for a precise statement. The center of the subalgebra $\U^{\circ}_q(\gl_n)$
is smaller than that of $\U_q(\gl_n)$; it turns out to be isomorphic to the algebra
of symmetric polynomials $\CC[q^{2\ell_1},\dots,q^{2\ell_n}]^{\Sym_n}$;
cf. \cite[Thm~2]{b:cq}, \cite[Thm~1.3]{jw:cr}.

A quantum version of the representation of $\gl_n$ in differential operators
for the reflection equation algebra was given by
Gurevich, Pyatov and Saponov~\cite{gps:bw}.
Elements of $\U^{\circ}_q(\gl_n)$ are represented by those of the {\em braided Weyl algebra}
and the corresponding $q$-analogues of the Capelli
identities associated with row and column Young diagrams were produced in \cite{gps:mc}.

We will rely on the Schur--Weyl duality between $\U_q(\gl_n)$ and the Hecke algebra
$\Hc_m$ to introduce the {\em $q$-immanant polynomials} $\SSb_{\mu}(z)$
by the equivalent formulas \eqref{imm} and \eqref{immsa} below, where they are written in terms of
a standard $\mu$-tableau $\Uc$, but they only depend on $\mu$ as shown in Theorem~\ref{thm:qimman}.
The coefficients of the polynomials $\SSb_{\mu}(z)$
are contained in the subalgebra $\U^{\circ}_q(\gl_n)$ and they are central in $\U_q(\gl_n)$.
We show that for any fixed value of $z$, as $\mu$ runs over diagrams with at most $n$ rows,
the elements $\SSb_{\mu}(z)$
form a basis of the center of $\U^{\circ}_q(\gl_n)$, and
their Harish-Chandra images
coincide with the factorial Schur polynomials
(Theorem~\ref{thm:qimman}).

The Harish-Chandra images of the constant terms $\SSb_{\mu}(0)$ are
the classical Schur polynomials, thus implying that the constant terms coincide with
the central elements provided by a general construction of Drinfeld~\cite{d:ac}
and Reshetikhin~\cite{r:qh}
based on the universal $R$-matrix; see also Etingof~\cite[Thm~1]{e:ce}.

We will set $\SSb_{\mu}=\SSb_{\mu}(z)$ for $z=(q-\qin)^{-1}$ and call these elements
the $q$-{\em immanants}. This term is justified by the fact that the limit value
of $\SSb_{\mu}$ as $q\to 1$ coincides with Okounkov's quantum immanant $\SSb_{\mu}$ for $\gl_n$
\cite{o:qi}. Moreover,
the images of the $q$-immanants $\SSb_{\mu}$
in the braided Weyl algebra take a factorized form thus yielding $q$-analogues of the
higher Capelli identities (Theorem~\ref{thm:capelli} and Corollary~\ref{cor:imimman}).
In the limit $q\to 1$ we recover the identities for $\gl_n$ along with
the Harish-Chandra images of the quantum immanants given in \cite{o:qi} and \cite{oo:ss}; see
Corollary~\ref{cor:hci}.

As a particular case of the $q$-immanant polynomials associated with the column diagrams we recover
the central elements of $\U_q(\gl_n)$ whose Harish-Chandra images are
the elementary symmetric polynomials. This allows us
to relate
them to the family of central elements $\tr_q\tss L^m$ of $\U_q(\gl_n)$
by an application of the quantum Liouville formula
of \cite{br:nb} and \cite{jlm:eq},
and to reproduce the Newton identities (Theorem~\ref{thm:newton})
previously given in \cite{gps:hs} and \cite{jw:cr}.

\section{Notation and definitions}
\label{sec:dch}

We regard $q$ as a nonzero complex number which is not
a root of unity.
As we defined in the Introduction, the {\em quantized enveloping algebra\/} $\U_q(\gl_n)$
is associated with the $R$-matrix
\beql{R}
R=q\ts\sum_{i}e_{ii}\ot e_{ii}+\sum_{i\ne j}e_{ii}\ot e_{jj}
+(q-\qin)\sum_{i< j}e_{ij}\ot e_{ji}\in \End\CC^n\ot \End\CC^n,
\eeq
where the $e_{ij}\in \End\CC^n$ denote the standard matrix units, while
$L^+$ and $L^-$ are the matrices
\ben
L^{\pm}=\sum_{i,j}e_{ij}\ot l^{\pm}_{ij}\in\End\CC^n\ot \U_q(\gl_n).
\een
We use a subscript to indicate a copy of the matrix in the multiple
tensor product algebras of the form
\beql{multtpr}
\underbrace{\End\CC^n\ot\dots\ot\End\CC^n}_m\ot\U_q(\gl_n),
\eeq
so that
\ben
L^{\pm}_a=\sum_{i,j=1}^n 1^{\ot (a-1)}\ot e_{ij}\ot 1^{\ot (m-a)}\ot l^{\pm}_{ij}.
\een
Furthermore, for an element
\ben
C=\sum_{i,j,r,s=1}^n c^{}_{ijrs}\ts e_{ij}\ot e_{rs}\in
\End \CC^n\ot\End \CC^n
\een
and any two indices $a,b\in\{1,\dots,m\}$ such that $a\ne b$,
we denote by $C_{a\tss b}$ the element of the algebra $(\End\CC^n)^{\ot m}$ with $m\geqslant 2$
given by
\beql{cab}
C_{a\tss b}=\sum_{i,j,r,s=1}^n c^{}_{ijrs}\ts (e_{ij})_a\tss (e_{rs})_b,
\qquad
(e_{ij})_a=1^{\ot(a-1)}\ot e_{ij}\ot 1^{\ot(m-a)}.
\eeq
For any $a\in\{1,\dots,m\}$ we will denote by $t_a$ the
partial transposition on the algebra \eqref{multtpr} which acts as the standard
matrix transposition on the
$a$-th copy of $\End \CC^n$ and as the identity map on all the other tensor factors.
The permutation operator $P$ is defined by
\beql{p}
P=\sum_{i,j=1}^{n}e_{ij}\ot e_{ji}\in
\End \CC^n\ot\End \CC^n.
\eeq
Given a matrix $X$, we will use
the $R$-matrix $\cR=P\tss R$ and follow \cite[Sec.~3.2]{i:qg}
to introduce another subscript notation
for elements of the algebra \eqref{multtpr} by setting
$X_{\ol 1}=X_{1}$ and
\ben
X_{\ol {k}}=\cR_{k-1}\dots \cR_{1}X_1\cR_{1}^{-1}\dots \cR_{k-1}^{-1},
\qquad k\geqslant 2,
\een
where we use the abbreviation $\cR_{i}=\cR_{i\ts i+1}$. Note that $X_{\ol {k}}=X_{k}$
in the specialization $q=1$.

Let $\mu$ be a Young diagram with $m$ boxes. If $\alpha=(i, j)$
is a box of $\mu$, then its {\em content} is $c(\al)=j-i$.
A tableau $\Uc$
of shape $\mu\vdash m$
(or a $\mu$-tableau)
is obtained by
filling in the boxes of the diagram with positive integers.
The entry in the box $\al\in\mu$ will be denoted by $\Uc(\al)$.
The tableau is called {\em semistandard}
if the entries weakly increase along each row from left to right
and strictly increase in each column from top to bottom.
We write $\sh(\Uc)=\mu$ if the shape of $\Uc$ is $\mu$.
A $\mu$-tableau with entries in
$\{1,\dots,m\}$ which are filled in the boxes bijectively
is called {\em standard}
if its entries strictly increase along the rows and down the columns.
We will denote by $f_{\mu}$ the number of standard tableaux of shape $\mu$.
Given a standard $\mu$-tableau $\Uc$, we let $c_{k}(\Uc)$ denote the content $j-i$
of the box $(i,j)$ of $\mu$ occupied by $k$ in $\Uc$.

We let $\Hc_m$ denote
the {\em Hecke algebra} over $\CC$ which is generated by elements
$T_1,\dots,T_{m-1}$ subject to the relations
\ben
\bal
(T_i-q)(T_i+q^{-1})=0,\\
T_{i}T_{i+1}T_{i}=T_{i+1}T_{i}T_{i+1},\\
T_iT_j=T_jT_i \quad\text{for\ \  $|i-j|>1$}.
\eal
\een
For $i=1, \dots, m-1$
we let $\si_{i}=(i, i+1)$ be the adjacent transpositions in the symmetric group $\Sym_{m}$.
Choose a reduced decomposition
$\si=\si_{i_{1}} \dots \si_{i_{l}}$ of any element
$\si \in \Sym_{m}$ and set $T_{\si}=T_{i_{1}} \dots T_{i_{l}}$.
The elements $T_{\si}$ parameterized by $\si\in\Sym_m$ form a basis of $\Hc_{m}$.

The irreducible representations of $\Hc_{m}$ over $\CC$ are parameterized by partitions $\mu$
of $m$.
The corresponding representation has a basis parameterized by the standard
$\mu$-tableaux.
The Hecke algebra $\Hc_m$ is semisimple and isomorphic to the direct sum of matrix algebras
\beql{isomhecke}
\Hc_m \cong \bigoplus_{\mu \vdash m} \Mat_{f_{\mu}}(\CC).
\eeq
The diagonal matrix units $e_{\Uc}=e_{\Uc \Uc}\in \Mat_{f_{\mu}}(\CC)$
with $\sh(\Uc)=\mu$ are {\em primitive idempotents} of $\Hc_m$.
They can be expressed explicitly in terms
of the pairwise commuting {\em Jucys--Murphy elements} $y_1,\dots,y_m$ of $\Hc_{m}$
or in terms of the {\em fusion procedure}, as reviewed e.g. in \cite{jlm:qs}.
By definition, $y_1=1$ and
\ben
y_k=1+(q-\qin)\ts
\big(T_{(1,k)}+T_{(2,k)}+\dots+T_{(k-1,k)}\big),\qquad k=2,\dots,m,
\een
where the elements $T_{(i,j)}\in\Hc_m$ are
associated with the transpositions $(i,j)\in\Sym_m$; see
\cite{c:ni}, \cite{dj:bi}. Equivalently, the Jucys--Murphy elements
can be defined by
\beql{jmind}
y_{k+1} = T_k\dots T_2\tss T_1^2\tss T_2\dots T_k
\qquad\text{for}\quad k=1,\dots,m-1.
\eeq
We have the relations
\beql{jmev}
y_k \ts e_{\Uc} =  e_{\Uc} \ts y_k = q^{2\tss c_k(\Uc)}\ts e_{\Uc}
\qquad\text{for}\quad k=1,\dots,m.
\eeq

\section{The $q$-immanants and eigenvalues}
\label{sec:qie}

The Hecke algebra $\Hc_{m}$ acts on the tensor product space
$(\CC^{n})^{\ot m}$
by
\beql{haact}
T_{k}\mapsto \cR_{k},\qquad k=1,\dots,m-1.
\eeq
We let $\cR_{\si}$ denote the image of the basis element $T_{\si}$.
Given a standard tableau $\Uc$ of shape $\mu\vdash m$,
we will let $\Ec_{\Uc}$ denote the image of the primitive idempotent $e_{\Uc}$
under the action \eqref{haact},
and working with the algebra \eqref{multtpr}
introduce the polynomial $\SSb_{\Uc}(z)$ in $z$
with coefficients in $\U_q(\gl_n)$
by
\begin{multline}\label{imm}
\SSb_{\Uc}(z)=\tr^{}_{1,\dots,m}\ts (L^+_1+zq^{-2\tss c_1(\Uc)}L^-_1)
\dots (L^+_m+zq^{-2\tss c_m(\Uc)}L^-_m)\\[0.3em]
{}\times(L^-_m)^{-1}\dots (L^-_1)^{-1}D_1\dots D_m \Ec_{\Uc},
\end{multline}
where\footnote{The matrix $D$ we use here differs from the one in \cite{jlm:qs} and \cite{jlm:eq}
by the scalar factor $q^{n-1}$.\label{footd}}
\beql{d}
D=\diag\big[1, q^{-2},\dots, q^{-2n+2}\big].
\eeq
In what follows, given a square matrix $X$ and a scalar $z$, expressions of the form $X+z$ are understood
as $X+z\tss 1$, where $1$ is the identity matrix of the same size as $X$.
Using the matrix $L=L^+(L^-)^{-1}$ and the $q$-{\em trace} $\tr_q X=\tr\ts DX$, we can write
$\SSb_{\Uc}(z)$ in the following equivalent form.

\bpr\label{prop:equiqm}
We have
\beql{immsa}
\SSb_{\Uc}(z)=\tr^{}_{q\ts (1,\dots,m)}\ts (L_{\ol 1}+zq^{-2\tss c_1(\Uc)})\dots
(L_{\ol m}+zq^{-2\tss c_m(\Uc)})\ts \Ec_{\Uc},
\eeq
assuming that the $q$-trace is taken over the copies of $\End\CC^n$ labelled by $1,\dots,m$.
\epr

\bpf
The equivalence clearly follows from the identity
\beql{idol}
(L^+_1+z_1\tss L^-_1)
\dots (L^+_m+z_m\tss L^-_m)
\ts(L^-_m)^{-1}\dots (L^-_1)^{-1}=(L_{\ol 1}+z_1)\dots
(L_{\ol m}+z_m),
\eeq
which is valid for arbitrary parameters $z_1,\dots,z_m$. To verify \eqref{idol}
use the reverse induction on $k$ starting from the trivial case $k=m$ to show that
\begin{multline}
(L^+_k+z_k\tss L^-_k)
\dots (L^+_m+z_m\tss L^-_m)
\ts(L^-_m)^{-1}\dots (L^-_k)^{-1}\\[0.3em]
=(L_{k}+z_k)\cR_k (L_{k}+z_{k+1})\cR_k^{-1}\dots
\cR_{m-1}\dots \cR_k(L_k+z_m)\cR_k^{-1}\dots \cR_{m-1}^{-1}.
\non
\end{multline}
Assuming this relation for a given $k\geqslant 2$, write
\ben
L^+_{k-1}+z_{k-1}\tss L^-_{k-1}=(L_{k-1}+z_{k-1})L^-_{k-1}.
\een
The induction step will be completed by observing that
$L^-_{k-1}$ commutes with $\cR_p$ with $p\geqslant k$ and
that the following relation holds:
\beql{lot}
L^-_{k-1}\tss L_{k}=\cR_{k-1}L_{k-1}\cR_{k-1}^{-1}\tss L^-_{k-1}.
\eeq
It suffices to take $k=2$ to verify \eqref{lot}, and it follows by writing
\ben
\bal
L^-_{1}\tss L_{2}=L^-_{1}\tss L^+_{2}(L^-_{2})^{-1}&=\cR\tss L^+_{1}\tss L^-_{2} \cR^{-1}(L^-_{2})^{-1}
=\cR\tss L^+_{1}\tss (L^-_{1})^{-1}\tss L^-_{1}\tss L^-_{2} \cR^{-1}(L^-_{2})^{-1}\\[0.4em]
&=\cR\tss L^+_{1}\tss (L^-_{1})^{-1}\tss \cR^{-1} L^-_{1}\tss L^-_{2}(L^-_{2})^{-1}
=\cR\tss L_1\tss \cR^{-1} L^-_{1},
\eal
\een
where we used the consequences
\beql{cpll}
\cR\tss L^{\pm}_{1}\tss L^{\pm}_{2}=L^{\pm}_{1}\tss L^{\pm}_{2}\cR\Fand
\cR\tss L^{+}_{1}\tss L^{-}_{2}=L^{-}_{1}\tss L^{+}_{2}\cR
\eeq
of the defining relations \eqref{rtt}.
\epf

To state the main theorem, assume that the Young diagram $\mu$
has at most $n$ rows and recall from \cite{gg:nt} and \cite{m:sft}
that the {\em factorial Schur polynomial}
in variables $x_1,\dots,x_n$
associated with $\mu$ and the sequence of parameters $a=(a_1,a_2,\dots)$ is given by the formula
\ben
s_{\mu}(x_1,\dots,x_n\tss|\tss a)=\sum_{\text{sh}(\Tc)=\mu}\ts
\prod_{\al\in\mu}
\big(x_{\Tc(\alpha)}+a_{\Tc(\alpha)+c(\alpha)}\big),
\een
summed over semistandard tableaux $\Tc$ of shape $\mu$ with
entries in the set $\{1,\dots,n\}$.

For an $n$-tuple of integers $\la=(\la_1,\dots,\la_n)$,
the irreducible representation
$L_q(\la)$ of $\U_q(\gl_n)$
is generated by a nonzero vector $\xi$ such that
\begin{alignat}{2}
l^-_{ij}\ts\xi&=0 \quad &&\text{for} \quad
i>j,
\non\\
l^+_{ii}\ts\xi&=q^{\la_i}\tss\xi \quad
&&\text{for} \quad i=1,\dots, n.
\non
\end{alignat}
Any element $w$ of the center of the algebra $\U_q(\gl_n)$ acts by the scalar multiplication
in the module $L_q(\la)$. The corresponding eigenvalue $\chi(w)$, regarded as a function of the
components of $\la$, coincides with the Harish-Chandra image of $w$.

\bth\label{thm:qimman}
\begin{enumerate}
\item
All coefficients of the polynomial $\SSb_{\Uc}(z)$ belong to the center of the
quantized enveloping algebra $\U_q(\gl_n)$.
\item
Given a diagram $\mu$,
the polynomial
$\SSb_{\Uc}(z)$ does not depend on the standard $\mu$-tableau $\Uc$.
\item
The eigenvalue of\ \  $\SSb_{\Uc}(z)$ in the module $L_q(\la)$ coincides
with the factorial Schur polynomial $s_{\mu}(q^{2\ell_1},\dots,q^{2\ell_n}\tss|\tss a)$,
where $\ell_i=\la_i-i+1$
and $a=(z, z\tss q^{-2},z\tss q^{-4},\dots)$.
\end{enumerate}
\eth

\bpf
All arguments are parallel to the proofs of the corresponding results of \cite{jlm:qs}
concerning the quantum Sugawara operators. Parts 1 and 2 follow by
modifying the proof of Theorem~3.1 therein in a straightforward way
since the structure of the key formula
\eqref{imm} is the same as the definition (3.14) of the quantum Sugawara operators in \cite{jlm:qs}.
We will indicate the necessary modifications below in the notation used
in the present paper. For a standard tableau $\Uc$ of shape $\mu\vdash m$ we now set
\ben
L^{+}_{\Uc}(z)=(L^+_1+zq^{-2\tss c_1(\Uc)}L^-_1)
\dots (L^+_m+zq^{-2\tss c_m(\Uc)}L^-_m)\Fand
L^{-}_{\Uc}=L^-_1\dots L^-_m
\een
so that
$\SSb_{\Uc}(z)$ can be written as
\ben
\SSb_{\Uc}(z)=\tr^{}_{1,\dots,m}\ts L_{\Uc}(z)D_1\dots D_m \Ec_{\Uc}
\een
with $L_{\Uc}(z)=L^{+}_{\Uc}(z)(L^{-}_{\Uc})^{-1}$.

\ble\label{lem:LE}
We have the relations
\beql{leee}
L^{+}_{\Uc}(z)\tss\Ec_{\Uc}=\Ec_{\Uc}\tss L^{+}_{\Uc}(z)\tss\Ec_{\Uc},
\qquad
(L^{-}_{\Uc})^{-1}\Ec_{\Uc}=\Ec_{\Uc}\tss (L^{-}_{\Uc})^{-1}\Ec_{\Uc},
\eeq
and hence
\ben
L_{\Uc}(z)\tss \Ec_{\Uc}=\Ec_{\Uc}\tss L_{\Uc}(z)\tss\Ec_{\Uc}.
\een
\ele

\bpf
The coefficients of the series $L^+(z)$ used in the proof of \cite[Lemma~3.2]{jlm:qs}
generate a subalgebra of the quantum affine algebra known as the $q$-{\em Yangian}.
Since the matrix $L^++zL^-$ coincides with the image of
$L^+(z)$ under
the evaluation homomorphism from the $q$-Yangian to $\U_q(\gl_n)$, the first relation
in \eqref{leee} follows. Furthermore, by the first relation in \eqref{cpll},
$L^{-}_{\Uc}$ commutes with the image of the Hecke algebra $\Hc_m$ under the action
\eqref{haact}, thus implying the second relation in \eqref{leee}.
\epf

Observe that by evaluating the $R$-matrix $R(x)$ used in \cite{jlm:qs} at $x=0$, we get
$q^{-1}R$ with $R$ defined in \eqref{R}. Hence, by taking $z=0$ in \cite[Lemma~3.3]{jlm:qs}
we come to the following lemma, where we use an additional tensor factor $\End\CC^n$
in \eqref{multtpr} labelled by $0$.

\ble\label{lem:RE}
We have the relations
\ben
\Ec_{\Uc}R_{0m}\dots
R_{01}\Ec_{\Uc}=
R_{0m}\dots R_{01}\Ec_{\Uc}
\een
and
\ben
\Ec_{\Uc}R_{01}^{-1}\dots
R_{0m}^{-1}\Ec_{\Uc}=
R_{01}^{-1}\dots
R_{0m}^{-1}\Ec_{\Uc}.
\vspace{-0.6cm}
\een
\qed
\ele

To show that
the coefficients of the polynomial $\SSb_{\Uc}(z)$ belong to the center of $\U_q(\gl_n)$,
we will verify that it commutes with $L^{+}_0$ and $L^-_0$. By using \eqref{rtt}
we get
\ben
L^+_0 L_{\Uc}(z)D_1\dots D_m\Ec_{\Uc}=
R_{01}^{-1}\dots
R_{0m}^{-1}\ts L_{\Uc}(z)
R_{0m}\dots
R_{01}
D_1\dots D_m\Ec_{\Uc}L^+_0.
\een
To conclude that
\beql{los}
L^+_0\tss \SSb_{\Uc}(z)=\SSb_{\Uc}(z)\tss L^+_0
\eeq
we need
to show that the trace
\beql{trace}
\tr^{}_{1,\dots,m}\ts
R_{01}^{-1}\dots
R_{0m}^{-1}\ts L_{\Uc}(z)
\tss R_{0m}\dots
R_{01}
D_1\dots D_m\Ec_{\Uc}
\eeq
equals $\SSb_{\Uc}(z)$. The product $D_1\dots D_m$ commutes with both the elements
of the symmetric group $\Sym_m$ and Hecke algebra $\Hc_m$ acting on the tensor product
space $(\CC^n)^{\ot m}$. Hence, using the relation $\Ec_{\Uc}^2=\Ec_{\Uc}$ and
applying Lemmas~\ref{lem:LE} and \ref{lem:RE}, we bring \eqref{trace} to the form
\ben
\tr^{}_{1,\dots,m}\ts
R_{01}^{-1}\dots
R_{0m}^{-1}
\Ec_{\Uc}L_{\Uc}(z)\Ec_{\Uc} R_{0m}\dots
R_{01}\tss\Ec_{\Uc}D_1\dots D_m.
\een
Set
\ben
X= R_{01}^{-1}\dots
R_{0m}^{-1}
\Ec_{\Uc}L_{\Uc}(z)\Ec_{\Uc},
\qquad
Y=  R_{0m}\dots
R_{01}\tss\Ec_{\Uc}D_1\dots D_m
\een
and note the identity
\ben
\tr^{}_{1,\dots,m}\ts XY=\tr^{}_{1,\dots,m}\ts X^{t_1\dots t_m}Y^{t_1\dots t_m}.
\een
Using the relation obtained by
applying the transposition $t_1\dots t_m$ to both sides of the second relation in Lemma~\ref{lem:RE}
we get
\begin{multline}
\tr^{}_{1,\dots,m}\ts X^{t_1\dots t_m}Y^{t_1\dots t_m}\\
{}=\tr^{}_{1,\dots,m}\ts\Ec_{\La}^{t_1\dots t_m}L_{\La}^{t_1\dots t_m}(w)
\big(R_{01}^{-1}\big)^{t_1}\dots
\big(R_{0m}^{-1}\big)^{t_m}\tss D_1\dots D_m\tss R_{0m}^{\ts t_m}\dots
R_{01}^{\ts t_1}.
\non
\end{multline}
By setting $x=0$ in the crossing symmetry relations \cite[(3.8)]{jlm:qs}, we obtain
the identities
\beql{csfi}
\big(R_{0k}^{-1}\big)^{t_k} D_k R_{0k}^{\ts t_k}=D_k
\Fand
R_{0k}^{\ts t_k}\tss D_k\big(R_{0k}^{-1}\big)^{t_k}=D_k
\eeq
for $k=1,\dots,m$ which are also easy to verify directly.
By using the first of these identities, we conclude that
$
\tr^{}_{1,\dots,m}\ts X^{t_1\dots t_m}Y^{t_1\dots t_m}
$
coincides with $\SSb_{\Uc}(z)$ thus completing the verification of \eqref{los}.
To verify its counterpart for $L^-_0$, note the equivalent form
of the second relation in \eqref{rtt}:
\ben
\wt R\ts L^-_1L^+_2=L^+_2 L^-_1\wt R,
\een
where $\wt R=(R_{21})^{-1}$. Furthermore, the $R$-matrix $R(x)$ used in \cite{jlm:qs}
can be written as
\ben
R(x)=\frac{f(x)}{q-\qin x}\ts (R-x\ts \wt R);
\een
see \cite[Sec.~3]{jlm:qs}. Hence, the relations of Lemma~\ref{lem:RE} hold in the same form
for the $R$-matrix $\wt R$ instead of $R$, as implied by cancelling numerical factors
in the relations of \cite[Lemma~3.3]{jlm:qs} and setting $w=0$.
Thus the above calculation proving \eqref{los} can be applied to $L^-_0$ instead of $L^+_0$
which only requires the replacement of the $R$-matrix $R$ with $\wt R$. In the last step
we apply the second identity in \eqref{csfi}. This completes
the proof of Part~1 of Theorem~\ref{thm:qimman}.

To prove Part~2 of Theorem~\ref{thm:qimman}, it will be enough to verify the relation
\beql{lalapr}
\tr^{}_{1,\dots,m}\ts L_{\Uc}(z)D_1\dots D_m \Ec_{\Uc}
=\tr^{}_{1,\dots,m}\ts L_{\Uc'}(z)D_1\dots D_m \Ec_{\Uc'},
\eeq
where $\Uc'=\si_k\ts\Uc$ is also a standard tableau for $k\in\{1,\dots,m-1\}$.
The calculation verifying the corresponding relation in the proof of Theorem~3.1 in
\cite{jlm:qs} relies on the key formula
\beql{lrch}
L_{\Uc}(z)\cR_k(q^{-2d_k}) = \cR_k(q^{-2d_k})L_{\Uc'}(z),
\eeq
where $d_k=c_{k+1}(\Uc)-c_k(\Uc)$. This formula holds for our settings as well,
as implied by \eqref{cpll} and the observation we made above that
the matrix $L^++zL^-$ coincides with the image of the series
$L^+(z)$ under
the evaluation homomorphism from the $q$-Yangian to $\U_q(\gl_n)$.
No changes are needed for
the rest of the calculations in the proof of Theorem~3.1 in
\cite{jlm:qs}.

Finally, to prove Part 3, we will rely on the arguments given
in the proof of \cite[Theorem~4.2]{jlm:qs}. We will outline the key steps and
indicate the necessary changes.

By Part~2, we can set $\SSb_{\mu}(z)=\SSb_{\Uc}(z)$. We will also use the
formula
\beql{sla}
\SSb_{\mu}(z)=\tr^{}_{1,\dots,m}\ts T_{\mu}(z),\qquad
T_{\mu}(z)=\frac{1}{f_{\mu}}\sum_{\sh(\Uc)=\mu}\ts L_{\Uc}(z)D_1\dots D_m\Ec_{\Uc},
\eeq
summed over standard $\mu$-tableaux $\Uc$. The following lemma
is proved in our settings in the same way as \cite[Lemma~4.3]{jlm:qs}.

\ble\label{lem:RT}
For any $\si\in \Sym_m$ and $\mu \vdash m$ we have
\ben
\cR_{\si}{T}_{\mu}(z)={T}_{\mu}(z)\cR_{\si}.
\vspace{-0.7cm}
\een
\qed
\ele

The subsequent calculations in the proof of Theorem~4.2 in \cite{jlm:qs} require
no essential changes, since the
highest weight conditions of the module $L_q(\la)$ are consistent with the ordering
of the generators of the quantum affine algebra in \cite[Sec.~4]{jlm:qs}, while
Lemmas~4.4--4.6 providing certain relations in the Hecke algebra $\Hc_m$ are used
in the same form. As a result, the eigenvalue of $\SSb_{\Uc}(z)$ in $L_q(\la)$
is given by the formula in \cite[Theorem~4.2]{jlm:qs}
which reads
\ben
\SSb_{\Uc}(z)\mapsto \sum_{\sh(\mathcal{T})=\mu}
\prod_{\al\in \mu}x^{}_{\mathcal{T}(\al)}(zq^{-2c(\al)}),
\een
summed over semistandard tableaux $\mathcal{T}$ of shape $\mu$ with entries in $\{1,2,\dots, n\}$,
where the series $x_i(z)$ now takes the form
\ben
x_i(z)=q^{n-2i+1}\ts (q^{2\la_i}+z)
\een
by replacing $l^{+}_{i}(z):=q^{\la_i}+z\tss q^{-\la_i}$ and $l^{-}_{i}(z):=q^{-\la_i}$.
The required formula for the eigenvalue then follows by taking into account the extra factor $q^{n-1}$
used in the definition of the matrix $D$; see footnote \ref{footd}.
\epf

By Part 2 of Theorem~\ref{thm:qimman}, we can define
the $q$-{\em immanant polynomial} associated with $\mu$ by
\beql{qimmsa}
\SSb_{\mu}(z)=\tr^{}_{q\ts (1,\dots,m)}\ts (L_{\ol 1}+zq^{-2\tss c_1(\Uc)})\dots
(L_{\ol m}+zq^{-2\tss c_m(\Uc)})\ts \Ec_{\Uc},
\eeq
for any standard tableau $\Uc$ of shape $\mu\vdash m$.

\bco\label{cor:basis}
\begin{enumerate}
\item
The elements
\beql{qimmsaze}
\SSb_{\mu}(0)=\tr^{}_{q\ts (1,\dots,m)}\ts L_{\ol 1}\dots
L_{\ol m}\ts \Ec_{\Uc},
\eeq
with $\mu$ running over the Young diagrams with at most $n$ rows form a basis of the center of
the algebra $\U^{\circ}_q(\gl_n)$.
\item
The eigenvalue of the central element $\SSb_{\mu}(0)$ in the module $L_q(\la)$ coincides
with the Schur polynomial $s_{\mu}(q^{2\ell_1},\dots,q^{2\ell_n})$.
\end{enumerate}
\eco

\bpf
In the case $z=0$, the sequence of parameters $a$ consists of zeros so that Part 3
of Theorem~\ref{thm:qimman}
implies that the eigenvalue of $\SSb_{\mu}(0)$ in the module $L_q(\la)$ coincides
with the Schur polynomial $s_{\mu}(q^{2\ell_1},\dots,q^{2\ell_n})$ yielding Part~2.
Since the Schur polynomials form a basis
of the algebra of symmetric polynomials in $q^{2\ell_1},\dots,q^{2\ell_n}$, Part 1 also follows.
\epf

Corollary~\ref{cor:basis} implies that the elements \eqref{qimmsaze} coincide with
the central elements constructed in~\cite{d:ac}
and~\cite{r:qh}; see also \cite[Thm~1]{e:ce}. They also
belong to a more general family of
central elements of the reflection equation algebra constructed in \cite[Prop.~5]{i:qg}.

\section{Higher quantum Capelli identities}
\label{sec:hqc}

The higher Capelli identities for $\gl_n$ can be proved by
a simple inductive argument as given in \cite[Sec.~7.4]{m:yc} and
a similar method was used in \cite{gps:mc} to prove $q$-analogues
of those identities associated with the row and column diagrams. Here we will extend the results of
\cite{gps:mc}
to get the identities associated with arbitrary Young diagrams
involving the $q$-immanant polynomials $\SSb_{\mu}(z)$.

\subsection{Action in the braided Weyl algebra}
\label{subsec:ab}

We will follow \cite{gps:bw} to
introduce the {\em braided Weyl algebra} $\Wc_n$ generated by two families
of elements $m_{ij}$ and $\di_{ij}$ with $1\leqslant i,j\leqslant n$ subject to the defining relations
\begin{align}
\cR M_1\cR M_1&=M_1\cR M_1\cR,
\label{mm}\\
\cR^{-1}\Dc_1\cR^{-1}\Dc_1&=\Dc_1\cR^{-1}\Dc_1\cR^{-1},
\label{dd}\\
\Dc_1\cR M_1\cR&=\cR M_1\cR^{-1}\Dc_1+\cR,
\label{dm}
\end{align}
where $M=[m_{ij}]$ and $\Dc=[\di_{ij}]$. We will see below in Sec.~\ref{subsec:hcq} that
at $q=1$ the algebra $\Wc_n$ coincides with
the algebra of polynomial coefficient differential operators
in the variables $m_{ij}$.

Introducing the matrix $K$ by
\beql{kl}
K=L+\frac{1}{q-\qin},
\eeq
we can rewrite the defining relations \eqref{re} of $\U^{\circ}_q(\gl_n)$ in the form
\beql{mre}
\cR\tss K_1\cR\tss K_1-K_1\cR\tss K_1\cR=\cR K_1-K_1\cR,
\eeq
known as the {\em modified reflection equation}; see e.g. \cite{gps:bw} and \cite[Sec.~3.2]{i:qg}.
According to \cite{gps:bw}, the mapping
\beql{diffre}
K\mapsto M\Dc
\eeq
defines a homomorphism $\U^{\circ}_q(\gl_n)\to \Wc_n$.

\bth\label{thm:capelli}
For any standard $\mu$-tableau $\Uc$,
the image of the element
\beql{eugen}
\Big(K_{\ol 1}+\frac{q^{-2c_1(\Uc)}-1}{q-\qin}\Big)\dots
\Big(K_{\ol m}+\frac{q^{-2c_m(\Uc)}-1}{q-\qin}\Big)\ts \Ec_{\Uc},
\eeq
under the homomorphism \eqref{diffre} is given by
\ben
a_{\mu}\ts M_{\ol 1}\dots M_{\ol m}\ts\Dc_{\ol m}\dots \Dc_{\ol 1}\ts \Ec_{\Uc}
\qquad\text{with}\qquad
a_{\mu}=\prod_{\al\in\mu} q^{-2c(\al)}.
\een
\eth

\bpf
We argue by induction on $m$ and write
\ben
\Ec_{\Uc}=\ga\ts\Ec_\Vc\ts\Ec_{\Uc},
\een
where $\Vc$ is the standard tableau obtained from $\Uc$ by
removing the box occupied by $m$ and $\ga$ is a nonzero constant.
We let $\nu$ be the shape of $\Vc$. The element \eqref{eugen} can
be written as
\ben
\ga\Big(K_{\ol 1}
+\frac{q^{-2c_1(\Uc)}-1}{q-\qin}\Big)\dots\Big(K_{\ol{m-1}}+\frac{q^{-2c_{m-1}(\Uc)}-1}{q-\qin}\Big)
\ts\Ec_\Vc\ts
\Big(K_{\ol m}+\frac{q^{-2c_m(\Uc)}-1}{q-\qin}\Big)\ts \Ec_{\Uc}.
\een
By the induction hypothesis,
its image under the homomorphism \eqref{diffre} equals
\begin{multline}
\ga\tss a_{\nu}\ts M_{\tss\ol 1}\dots M_{\ts\ol {m-1}}\ts\Dc_{\ts\ol {m-1}}\dots
\Dc_{\tss\ol 1}\ts \Ec_{\Vc}\Big(K_{\ol m}+\frac{q^{-2c_m(\Uc)}}{q-\qin}\Big)\ts \Ec_{\Uc}\\
=a_{\nu}\ts M_{\ol 1}\dots M_{\ts\ol {m-1}}\ts\Dc_{\ts\ol {m-1}}\dots
\Dc_{\tss\ol 1}\ts \Big(K_{\ol m}+\frac{q^{-2c_m(\Uc)}-1}{q-\qin}\Big)\ts \Ec_{\Uc},
\non
\end{multline}
where we keep notation $K_{\ol m}$ for the image of this element under \eqref{diffre}.
As was already calculated in \cite[Sec.~3]{gps:mc},
\ben
\Dc_{\tss\ol k}\tss K_{\ol m}=K_{\ol m}\tss \Dc_{\tss\ol k}\ts a_k
+ \Dc_{\tss\ol k}\ts b_k,\qquad k=1,\dots,m-1,
\een
with
\ben
a_k=\cR_{m-1}\dots\cR_{k+1}\cR_k^{-2}\cR_{k+1}^{-1}\dots \cR_{m-1}^{-1}
\fand
b_k=\cR_{m-1}\dots\cR_{k+1}\cR_k^{-1}\cR_{k+1}^{-1}\dots \cR_{m-1}^{-1}.
\een
Hence
\begin{multline}
\Dc_{\ts\ol {m-1}}\dots\Dc_{\ts\ol 1}\tss K_{\ol m}=K_{\ol m}
\ts \Dc_{\ts\ol {m-1}}\dots\Dc_{\ts\ol 1}\ts
a_{m-1}\dots a_1\\[0.3em]
{}+\Dc_{\ts\ol {m-1}}\dots\Dc_{\ts\ol 1}\ts(b_1+b_2a_1+\dots+b_{m-1}a_{m-2}\dots a_1).
\label{codi}
\end{multline}
Now write $K_{\ol m}=M_{\tss\ol m}\Dc_{\tss\ol m}$
and observe that
\ben
a_{m-1}\dots a_1=\cR_{m-1}^{-1}\dots\cR_{2}^{-1}\cR_1^{-2}\cR_{2}^{-1}\dots \cR_{m-1}^{-1}
\een
which coincides with the image of the inverse Jucys--Murphy element $y_m$ in \eqref{jmind}
under
the action \eqref{haact} of the Hecke algebra. Therefore by \eqref{jmev}
\beql{amsi}
a_{m-1}\dots a_1\ts \Ec_{\Uc}=q^{-2c_m(\Uc)}\ts \Ec_{\Uc}.
\eeq
Furthermore, since
\ben
\cR_k^{-1}=\frac{1}{q-\qin}(1-\cR_k^{-2})
\een
we have
\ben
b_k=\frac{1}{q-\qin}(1-a_k).
\een
Therefore, the coefficient of $\Dc_{\ts\ol {m-1}}\dots\Dc_{\ts\ol 1}$ in
\eqref{codi} simplifies via a telescoping sum to
\ben
\frac{1}{q-\qin}(1-a_{m-1}\dots a_1).
\een
This implies that by \eqref{amsi} the coefficient of
the product $M_{\ts\ol 1}\dots M_{\ts\ol {m-1}}\ts\Dc_{\ts\ol {m-1}}\dots
\Dc_{\ts\ol 1}\ts \Ec_{\Uc}$ in the final expression vanishes.
The proof is completed by observing that $a_{\mu}=a_{\nu}\ts q^{-2c_m(\Uc)}$.
\epf

The particular cases $\mu=(m)$ and $\mu=(1^m)$ of Theorem~\ref{thm:capelli}
are contained in \cite[Thm~3]{gps:mc}.

Note that the $q$-immanant $\SSb_{\mu}=\SSb_{\mu}((q-\qin)^{-1})$ coincides with
the $q$-trace of the element on the left hand side of \eqref{eugen}.
Hence Theorem~\ref{thm:capelli} implies the following.

\bco\label{cor:imimman}
The image of the $q$-immanant $\SSb_{\mu}$ under
the homomorphism \eqref{diffre} is given by
\ben
\SSb_{\mu}\mapsto a_{\mu}\ts \tr^{}_{q\ts (1,\dots,m)}\ts M_{\ol 1}\dots M_{\ol m}\ts
\Dc_{\ol m}\dots \Dc_{\ol 1}\ts \Ec_{\Uc}.
\een
\eco

\subsection{To the higher Capelli identities for $\gl_n$}

\label{subsec:hcq}

Consider the `semiclassical limit' $q\to 1$
in Theorem~\ref{thm:capelli} and Corollary~\ref{cor:imimman}.
Note that at the evaluation $q=1$ we have $\cR=P$ so that
relation \eqref{mm} becomes $M_1\tss M_2=M_2\tss M_1$ meaning that the elements $m_{ij}$
pairwise commute. Similarly, relation \eqref{dd} manifests the same property for the $\di_{ij}$,
while the third relation takes the form $\Dc_1M_2=M_2\Dc_1+P_{12}$ which is equivalent to
\ben
\di_{ij}\tss m_{kl}=m_{kl}\tss \di_{ij}+\de_{kj}\de_{il}.
\een
Thus, $\Wc_n$ becomes
the algebra of polynomial coefficient differential operators
in the variables $m_{ij}$ with $\di_{ij}=\di/\di m_{ji}$.

By taking $q=1$ in \eqref{mre} and replacing $K\mapsto E$, we recover the defining relations
\ben
E_1E_2-E_2E_1=E_1P-PE_1
\een
of the
universal enveloping algebra $\U(\gl_n)$ written for the generator matrix
\ben
E=\sum_{i,j=1}^n e_{ij}\ot E_{ij}\in\End\CC^n\ot \U(\gl_n)
\een
with the standard notation $E_{ij}$ of the basis elements of the Lie algebra $\gl_n$.
Therefore, mapping \eqref{diffre}
is a $q$-deformation of the representation of $\gl_n$ in the
polynomial coefficient differential operators.

By introducing the entries of the matrix $K=[k_{ij}]$, for the action on the
highest vector $\xi\in L_q(\la)$ we get
\begin{alignat}{2}
k_{ij}\ts\xi&=0 \quad &&\text{for} \quad
i>j,
\non\\
k_{ii}\ts\xi&=\vk_i\ts\xi \quad
&&\text{for} \quad i=1,\dots, n\quad \text{with} \quad \vk_i=q^{2\tss\la_i}+\frac{1}{q-\qin}.
\non
\end{alignat}
For an arbitrary $n$-tuple of complex numbers $\vk=(\vk_1,\dots,\vk_n)$
we let $L(\vk)$ denote the irreducible highest
weight representation of $\gl_n$ generated by a nonzero vector $\xi$ such that
\begin{alignat}{2}
E_{ij}\ts\xi&=0 \quad &&\text{for} \quad
i>j,
\non\\
E_{ii}\ts\xi&=\vk_i\ts\xi \quad
&&\text{for} \quad i=1,\dots, n.
\non
\end{alignat}

By using the matrix $K$ defined in \eqref{kl} instead of $L$, we can write the
$q$-immanant in the form
\beql{keugen}
\SSb_{\mu}=\tr^{}_{q\ts (1,\dots,m)}\ts\Big(K_{\ol 1}+\frac{q^{-2c_1(\Uc)}-1}{q-\qin}\Big)\dots
\Big(K_{\ol m}+\frac{q^{-2c_m(\Uc)}-1}{q-\qin}\Big)\ts \Ec_{\Uc},
\eeq
which has a well-defined limit as $q\to 1$ equal to the quantum immanant
\beql{lism}
\SSb_{\mu}=\tr^{}_{1,\dots,m}\ts\big(E_1-c_1(\Uc)\big)\dots
\big(E_m-c_m(\Uc)\big)\ts \Ec_{\Uc},
\eeq
where $\Ec_{\Uc}$ now denotes the corresponding primitive idempotent for $\Sym_m$.
We thus recover the higher Capelli identities
for $\gl_n$ together with the Harish-Chandra images of the quantum $\mu$-immanants
\eqref{lism}
as given in \cite{o:qi}, \cite{oo:ss}; cf. \cite{jlz:gc} and \cite{w:so}.

\bco\label{cor:hci}
\begin{enumerate}
\item
The image of the element \eqref{lism}
under the action $E\mapsto M\Dc$ of $\gl_n$
in the polynomial coefficient differential operators
equals
\beql{mddi}
\tr^{}_{1,\dots,m}\ts M_{1}\dots M_{m}\ts
\Dc_{1}\dots \Dc_{m}\ts \Ec_{\Uc}.
\eeq
\item
The eigenvalue of $\SSb_{\mu}$ in the module $L(\vk)$
is given by
\ben
\sum_{\text{sh}(\Tc)=\mu}\ts
\prod_{\al\in\mu}
\big(\vk^{}_{\Tc(\alpha)}-c(\alpha)\big),
\een
summed over semistandard tableaux $\Tc$ of shape $\mu$ with
entries in the set $\{1,\dots,n\}$.
\end{enumerate}
\eco

\bpf
Part 1 is immediate from Theorem~\ref{thm:capelli}. Furthermore, by Part~3 of Theorem~\ref{thm:qimman}
with $z=(q-\qin)^{-1}$, the eigenvalue of $\SSb_{\mu}$ in $L_q(\la)$ is found by
\ben
\sum_{\text{sh}(\Tc)=\mu}\ts
\prod_{\al\in\mu}
\big(q^{2\tss\ell_{\Tc(\alpha)}}+\frac{q^{-2\tss \Tc(\alpha)-2\tss c(\al)+2}}{q-\qin}\big).
\een
Write this expression in terms of the parameters $\vk_i$ to get
\ben
\sum_{\text{sh}(\Tc)=\mu}\ts
\prod_{\al\in\mu}\ts q^{-2\tss \Tc(\alpha)+2}\ts
\Big(\vk^{}_{T(\al)}+\frac{q^{-2\tss c(\al)}-1}{q-\qin}\big).
\een
Part~2 now follows by taking the limit $q\to 1$.
\epf

\section{Newton identities}
\label{sec:ni}

As an application of Theorem~\ref{thm:qimman}, we will derive the Newton identities connecting
two families of central elements in the algebra $\U_q(\gl_n)$.
Our starting point is the quantum Liouville formula stated in \cite[Eq. (4.28)]{br:nb},
a proof is given in \cite{jlm:eq}. To recall the formula, consider
the {\em quantum loop algebra $\U_q(\wh\gl_n)$}
generated by elements
\ben
l^+_{ij}[-r],\qquad l^-_{ij}[r]\qquad\text{with}\quad 1\leqslant i,j\leqslant n,\qquad r=0,1,\dots,
\een
subject to the defining relations
\begin{align}
l^+_{ji}[0]&=l^-_{ij}[0]=0\qquad&&\text{for}\qquad 1\leqslant i<j\leqslant n,
\non\\
l^+_{ii}[0]\ts l^-_{ii}[0]&=l^-_{ii}[0]\ts l^+_{ii}[0]=1\qquad&&\text{for}\qquad i=1,\dots,n,
\non
\end{align}
and
\begin{align}
R(u/v)L_1^{\pm}(u)L_2^{\pm}(v)&=L_2^{\pm}(v)L_1^{\pm}(u)R(u/v),
\label{RLL}\\[0.3em]
R(u/v)L_1^{+}(u)L_2^{-}(v)&=L_2^{-}(v)L_1^{+}(u)R(u/v),
\label{RLLpm}
\end{align}
with
\beql{rx}
R(x)=R-\frac{x\tss(q-q^{-1})}{x-1}P,
\eeq
where $R$ and $P$ are defined in \eqref{R} and \eqref{p}.
In the defining relations we consider the matrices $L^{\pm}(u)=\big[\tss l^{\pm}_{ij}(u)\big]$,
whose entries are formal power series in $u$ and $u^{-1}$,
\ben
l^{+}_{ij}(u)=\sum_{r=0}^{\infty}l^{+}_{ij}[-r]\tss u^r,\qquad
l^{-}_{ij}(u)=\sum_{r=0}^{\infty}l^{-}_{ij}[r]\tss u^{-r}.
\een
The quantum Liouville formula reads
\beql{liouvpm}
z^{\pm}(u)=\frac{\qdet L^{\pm}(uq^{2})}{\qdet L^{\pm}(u)}
\eeq
for the series $z^{\pm}(u)$ defined by\footref{footd}
\ben
z^{\pm}(u)\ts\tr\ts D=\tr\ts DL^{\pm}(uq^{2n})L^{\pm}(u)^{-1},
\een
while $\qdet L^{\pm}(u)$ is the quantum determinant
\ben
\qdet L^{\pm}(u)=
\sum_{\si\in \Sym_n} (-q)^{-l(\si)} \ts l^{\pm}_{1\si(1)}(u)\cdots
l^{\pm}_{n\si(n)}(u q^{2n-2}),
\een
where $l(\si)$ denotes the length of the permutation $\si$. All coefficients of the series
$\qdet L^{\pm}(u)$ and $z^{\pm}(u)$ belong to the center of $\U_q(\wh\gl_n)$.

The quantized enveloping algebra $\U_q(\gl_n)$ can be regarded as a subalgebra of $\U_q(\wh\gl_n)$
via the embedding $L^{\pm}\mapsto L^{\pm}[0]$, where $L^{\pm}[0]=[l_{ij}^{\pm}[0]]$. We will
identify $L^{\pm}$ with $L^{\pm}[0]$.
Introduce the matrices $M^{\pm}(u)=[m_{ij}^{\pm}(u)]$ by
\ben
M^{\pm}(u)=L^{\pm}(u)(L^-)^{-1}
\een
and denote by $\U^{\circ}_q(\wh\gl_n)$ the subalgebra of $\U_q(\wh\gl_n)$
generated by all coefficients of the series $m_{ij}^{\pm}(u)$.
These matrices satisfy
the relations\footnote{Each relation in \eqref{braid} defines the ``braided Yangian''
in the terminology of \cite{gs:hs}. In our settings, this is just a subalgebra of the quantum loop
algebra.}
\begin{align}\label{braid}
\cR(u/v)M^{\pm}_1(u)\cR\tss M^{\pm}_1(v)&=M^{\pm}_1(v)\cR\tss M^{\pm}_1(u)\cR(u/v),\\[0.3em]
\label{braidpm}
\cR(u/v)M^{+}_1(u)\cR\tss M^{-}_1(v)&=M^{-}_1(v)\cR\tss
M^{+}_1(u)\cR(u/v),
\end{align}
where $\cR(x)=P\tss R(x)$. This follows easily with the use of \eqref{cpll} and the relations
\beql{rlul}
\cR\ts L_1^{\pm}(u)L_2^{-}=L_1^{-}L_2^{\pm}(u)\cR
\eeq
implied by \eqref{RLL} and \eqref{RLLpm}. The images of the matrices $M^{\pm}(u)$
under the evaluation homomorphism $\U_q(\wh\gl_n)\to \U_q(\gl_n)$
defined by
\beql{eval}
L^+(u)\mapsto L^+-L^- u,\qquad
L^-(u)\mapsto L^--L^+ u^{-1},
\eeq
are
given by
\beql{evalbr}
M^+(u)\mapsto L-u,\qquad
M^-(u)\mapsto 1-L u^{-1}.
\eeq

Returning to the Liouville formula \eqref{liouvpm}, observe
that it can be interpreted as a relation in the subalgebra $\U^{\circ}_q(\wh\gl_n)$
since the coefficients of both the series $z^{\pm}(u)$ and the ratio of the quantum determinants
can be regarded as elements of $\U^{\circ}_q(\wh\gl_n)$. This is clear for $z^{\pm}(u)$,
and we claim that
\beql{qdlm}
\qdet L^{\pm}(u)=\qdet M^{\pm}(u)\tss l_{11}^-\dots l_{nn}^-
\eeq
for series $\qdet M^{\pm}(u)$ with coefficients in $\U^{\circ}_q(\wh\gl_n)$. Indeed, we
have the relations defining the quantum determinant,
\ben
A^{(n)}_q L^{\pm}_1(uq^{2n-2})\cdots L^{\pm}_n(u)
=A^{(n)}_q\tss \qdet L^{\pm}(u),
\een
where $A^{(n)}_q$ is the $q$-antisymmetrizer; see e.g. \cite{jlm:eq}. In particular,
\beql{almm}
A^{(n)}_q L^{-}_1\cdots L^{-}_n
=A^{(n)}_q\tss l_{11}^-\dots l_{nn}^-,
\eeq
and the product $l_{11}^-\dots l_{nn}^-$ is central in $\U_q(\wh\gl_n)$.
On the other hand,
we have
\ben
L^{\pm}_1(uq^{2n-2})\cdots L^{\pm}_n(u)
=M^{\pm}_{\ol 1}(uq^{2n-2})\cdots M^{\pm}_{\ol n}(u)\ts L_1^-\dots L_n^-.
\een
This is verified in the same way as \eqref{idol} with the use
of \eqref{rlul} and the counterpart of \eqref{lot} given by
\ben
L^-_{k-1}\tss M^{\pm}_{k}(v)=\cR_{k-1}M^{\pm}_{k-1}(v)\cR_{k-1}^{-1}\tss L^-_{k-1}.
\een
Hence
\ben
A^{(n)}_q M^{\pm}_{\ol 1}(uq^{2n-2})\cdots M^{\pm}_{\ol n}(u)L_1^-\dots L_n^-=
A^{(n)}_q\tss \qdet L^{\pm}(u)
\een
so that \eqref{qdlm} follows by writing this relation in the form
\ben
A^{(n)}_q M^{\pm}_{\ol 1}(uq^{2n-2})\cdots M^{\pm}_{\ol n}(u)=
A^{(n)}_q\tss \qdet L^{\pm}(u)(L_n^-)^{-1}\dots (L_1^-)^{-1}
\een
and using \eqref{almm}.

By taking the corresponding version
\beql{liouvpmm}
z^{-}(u)=\frac{\qdet M^{-}(uq^{2})}{\qdet M^{-}(u)}
\eeq
of \eqref{liouvpm}, we find that the image of $z^{-}(u)$ under the evaluation map
in the second formula in \eqref{evalbr} is given by
\ben
1+(1-q^{-2})\ts\sum_{m=1}^{\infty}\ts \tr_q\tss L^m\tss u^{-m}.
\een
Performing the same calculation as in \cite[Sec.~3]{jlm:eq}
for the quantum determinant,
we conclude that for the eigenvalues of the central elements of $\U^{\circ}_q(\gl_n)$
in $L_q(\la)$ we have the identity
\beql{idliou}
1+(1-q^{-2})\ts\sum_{m=1}^{\infty}\ts \chi(\tr_q\tss L^m)\tss u^{-m}
=\frac{(1-q^{2\ell_1-2}u^{-1})\dots (1-q^{2\ell_n-2}u^{-1})}
{(1-q^{2\ell_1}u^{-1})\dots (1-q^{2\ell_n}u^{-1})},
\eeq
where, as before, $\ell_i=\la_i-i+1$ and $\chi(w)$ denotes the eigenvalue of a central element
$w$ in $L_q(\la)$.
By \eqref{qimmsaze}, the constant term $\SSb_{(1^m)}(0)$ associated with
the column diagram $(1^m)$ is given by
\ben
\SSb_{(1^m)}(0)=\tr^{}_{q\ts (1,\dots,m)}\ts L_{\ol 1}\dots L_{\ol m}\ts \Ec_{(1^m)}.
\een
Due to Part~2 of Corollary~\ref{cor:basis},
the eigenvalue of $\SSb_{(1^m)}(0)$ in $L_q(\la)$ is the elementary
symmetric polynomial
\ben
e_m(q^{2\ell_1},\dots,q^{2\ell_n})
=\sum_{1\leqslant i_1<\dots <i_m\leqslant n} q^{2(\ell_{i_1}+\dots+\ell_{i_m})}.
\een
Hence the denominator in \eqref{idliou}
is the eigenvalue of the generating function $E(u)$
of the central elements defined by
\ben
E(u)=\sum_{m=0}^n \tr_{q(1,\dots,m)}\tss
L_{\ol 1}\dots L_{\ol m}\ts \Ec_{(1^m)}\ts (-u)^{-m}.
\een
Thus, we have derived the following version of the Newton identities which
already appeared
in different forms in \cite{gps:hs} and \cite{jw:cr}.

\bth\label{thm:newton}
We have the identity
\ben
1+(1-q^{-2})\ts\sum_{m=1}^{\infty}\ts \tr_q\tss L^m\tss u^{-m}
=\frac{E(u\tss q^2)}{E(u)}.
\een
\eth


\bigskip
\bigskip

\small

\noindent
N.J.:\newline
Department of Mathematics\\
North Carolina State University, Raleigh, NC 27695, USA\\
jing@ncsu.edu

\vspace{5 mm}

\noindent
M.L.:\newline
School of Mathematical Sciences\\
South China University of Technology\\
Guangzhou, Guangdong 510640, China\\
mamliu@scut.edu.cn

\vspace{5 mm}

\noindent
A.M.:\newline
School of Mathematics and Statistics\newline
University of Sydney,
NSW 2006, Australia\newline
alexander.molev@sydney.edu.au

\end{document}